# THE DISCRETE AND THE CONTINUOUS: WHICH COMES FIRST?

T. N. Narasimhan
February 7, 2010

*In solving diffusion problems, it is common to consider the finite difference equation to be an approximation to the differential equation. Nevertheless, history shows that the finite difference equation is primitive and that the differential equation is its idealized representation designed to obtain solutions in algebraic form. The difference equation is logically consistent within itself, independent of the differential equation. The difference equation and the differential equation together constitute two powerful complementary tools, one providing numerical solutions to problems of arbitrary complexity on a case by case basis, and the other providing insights into classes of problems under idealized conditions.*

**Introduction**

The notion of a continuum and its description with a differential equation are ruling paradigms of mathematical physics. Over the past half a century, the digital computer has enabled numerical solution of complex physical problems defined over discrete domains. Commonly, such discretely defined problems are treated as approximations to continuum problems described by differential equations. The observational world is discrete and finite. If so, how may one assume that the discrete representation is an approximation to an idealized continuum and the associated differential equation? Which comes first, the discrete or the continuous? To find answers, it is necessary to examine the history of difference equations and differential equations. Accordingly, we examine eighteenth century investigations of probability, followed by the birth of the heat equation during early nineteenth century.

**Difference and Differential Equations in Probability**

Following the introduction of mathematical probability as a tool for describing outcomes of games of chance, the study of probability engaged the attention of distinguished mathematicians during the eighteenth century. Initially, these problems were investigated using combinatoric methods. For example, James Bernoulli[1] considered a problem involving two events, one of them having probability of occurrence p, and the other having probability of q = 1-p. He showed that the probability that the first event will occur m times and fail n times is equal to a certain term in the expansion of $(p+q)^\mu$, namely,

$$\frac{\mu!}{m!\,n!} p^m q^n, \qquad (1)$$

where $\mu = (m+n)$. Unfortunately, in the absence of computing devices, expressions such as (1) could not be numerically evaluated when $\mu$ was large.

To overcome this difficulty, they devised a strategy of evaluating these ratios using definite integrals. For this purpose, the technique of génératrice function was introduced during the middle of the 18$^{th}$ century. The central idea was to set up a power series expansion in which the required probability will form a coefficient. Having set up such an equation, the next step was to evaluate a desired coefficient by setting up an appropriate difference equation expressing the incremental change in the coefficient as $\mu$ is increased to $(\mu+1)$.

Pursuing this approach, Pierre Simon Laplace[2] considered a power series involving products of two variables t and t' and established the following finite difference expressions for the coefficient $y_{x,x'}$ of the product $t^x t'^{x'}$ :

$$\Delta^2 y_{x,x'} = y_{x+2,\,x'} - 2 y_{x+1,\,x'} + y_{x,x'}$$

$$\Delta' y_{x,x'} = y_{x,x'+1} - y_{x,x'} .$$

Note that x and x' are integers, and $\Delta x = \Delta x' = 1$. Therefore, the right hand sides of (2a) and (2b) denote typical finite difference expressions with $(\Delta x)^2$ and $\Delta x'$ implied as denominators.

Later, Laplace[3] used recursive reasoning to set up the partial difference equation,

$$\Delta^2 y_{x,x'} = \Delta' y_{x,x'} . \qquad (3)$$

Invoking infinitesimal reasoning, he then replaced the difference equation with the partial differential equation,

$$\frac{\partial^2 y}{\partial x^2} = \frac{\partial y}{\partial x'},\qquad(4)$$

where y is relative frequency (probability density) representing $y_{x,x'}$. Joseph Fourier[4] built on Laplace's solution to (4) and showed that that (4) is satisfied by normal distribution.

Clearly, the principal idea underlying differential equation (4) is a difference equation.

**Difference and Differential Equations in Heat Diffusion**

At the beginning of the nineteenth century, Fourier[5,6] embarked on setting up a partial differential equation describing transient heat flow in a solid. His rationale was to consider a small, discrete volume element in the interior of a solid, sum up heat fluxes flowing into and out of the bounding surfaces of the element over a small interval of time, and divide net accumulation of heat by the heat capacity of the element to represent average change in temperature of the element over the chosen time interval. To evaluate heat fluxes, he imagined that heat flowed between isothermal surfaces, with the quantity of heat per unit time being directly proportional to temperature difference and cross-sectional area, and inversely proportional to distance between two surfaces. Figure 1 (ref. 6, p. 97) illustrates his reasoning.

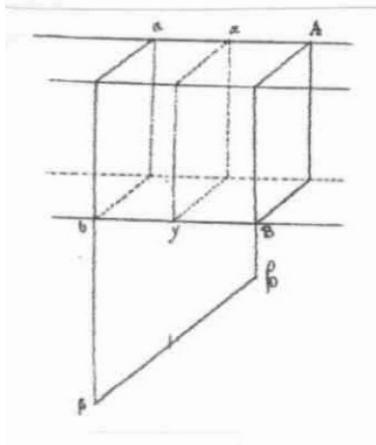

Figure 1:   Fourier's conceptualization of thermal conductivity using heat flow in a prism of unit cross sectional area (Fourier's symbols[6])

The isothermal surface on the left at x = α has temperature b, and that at x = A has temperature B. The intermediate isothermal surface at x has a temperature y. Assuming that under steady conditions, temperature varies linearly between the surfaces, Fourier expressed the quantity of heat F flowing across unit cross sectional area per unit time by,

$$F = K\frac{b-B}{\alpha-A} = -K\frac{dy}{dx},$$

where the abcissa is positive to the right. In (5), linear slope is treated as synonymous with gradient. Next, Fourier considered conservation of heat in a slab M of width Δx and cross section S, communicating with two neighboring slabs L and N (Figure 2), assuming linear variation of temperature over distances of the order of Δx.

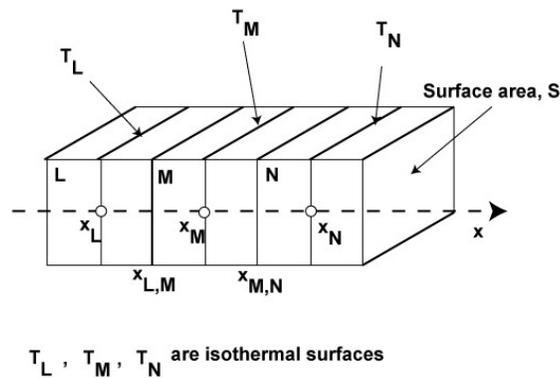

Figure 2: Conservation of heat in a slab M connected to two neighboring slabs in a prism

The average temperature $T_M$ of slab M at an instant of time is equal to the heat content of the slab divided by its heat capacity. Thus,

$$T_M = \frac{\int_{x_{L,M}}^{x_{M,N}} \rho c T(y) S dy}{S(x_{M,N} - x_{L,M})}.$$

But, since T varies linearly within the slab, the average temperature given by (6) occurs at the mid point of the slab, $x_M$. These considerations lead to the difference equation,

$$KS(\frac{T_L - T_M}{\Delta x} - \frac{T_M - T_N}{\Delta x}) = S \Delta x \rho c \frac{\Delta T_M}{\Delta t},$$

where the change in temperature $\Delta T_M$ denotes change in temperature of the isothermal surface passing through mid-point $x_M$ of the slab.

This finite difference equation describes one-dimensional transient heat flow in a prism in sufficient detail to obtain numerical solutions when boundary conditions and initial conditions are prescribed. Obviously, no computing devices capable of handling the numerical calculations involved in implementing (7) were available at the time of Fourier.

To convert (7) into a differential equation, Fourier divided through by volume, (S Δx), of the element, and set ($\Delta T_M /\Delta t$) = dT/dt for a sufficiently small Δt (ref. 6, p. 101-102) to arrive at,

$$K \frac{\partial^2 T}{\partial x^2} = \rho c \frac{\partial T}{\partial t}.$$

To extend (8) to bodies of more general shape such as a sphere or a cylinder, Fourier did not set up difference equations. Instead, he invoked infinitesimals in combination with (8), accounted for variation of surface area with x to formulate appropriate differential equations.

Fourier's Analytic Theory of Heat published in 1822 inspired Georg Simon Ohm, Adolf Fick, and James Clerk Maxwell to apply his heat flow model to electricity, molecular diffusion, and electromagnetism, each introducing novel conceptual advancements in the process. Ohm[7] proposed the concept of a resistance and expressed electrical flux in terms of a difference equation, rather than using a gradient as Fourier did. Fick[8,9] who conducted salt-diffusion experiments with a vessel in the shape of a truncated cone, used a generalized one-dimensional diffusion equation for flow in a tube with variable cross section,

$$K\left(\frac{d^2c}{dx^2} + \frac{1}{A}\frac{dA}{dx}\frac{dc}{dx}\right) = \frac{dc}{dt} ,$$

where c is concentration, and A is the area of cross section perpendicular to flow. The exact solution to this equation gives a generalized expression of Ohm's Law,

$$F = \frac{T_{in} - T_{out}}{R_{in,out}} = \frac{T_{in} - T_{out}}{\int_{x_{in}}^{x_{out}} \frac{dy}{KA(y)}} ,$$

where $R_{in,out}$ denotes the resistance between inlet and outlet. Maxwell[10,11], inspired by Michael Faraday's concept of lines of force around a magnet, visualized flow of an incompressible fluid in a resistive medium in terms of flow lines and a collection of flow tubes with orthogonal surfaces of equal pressure[12]. Under dynamical equilibrium, forces impelling fluid flow were assumed to be exactly balanced by resistive forces, and flux was effectively given by (10).

Together, the contributions of Ohm, Fick, and Maxwell enable a fully self-consistent finite-difference statement of the heat flow problem, independent of the differential equation. Recall that in setting up difference equation (7) for the prism, average temperature was associated with the isothermal surface exactly midway between interfaces bounding an element. However, when the area of cross section is variable, linear temperature variation along flow path cannnot be assumed, and average temperature cannot be associated with an isothermal surface midway between bounding interfaces. With this recognition, consider now three adjoining elements L, M, and N along a tube with variable cross section, as shown in Figure 3.

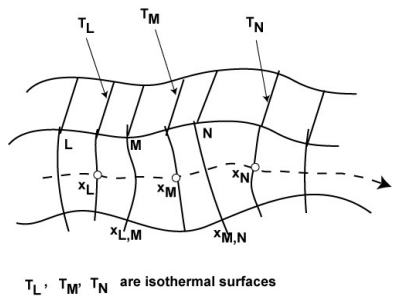

$T_L$, $T_M$, $T_N$ are isothermal surfaces

Figure 3: Flow tube of variable cross section with a curvilinear x axis

For element M bounded by interfaces at $x_{L,M}$ and $x_{M,N}$, average temperature at an instant t is given, analogous to (6) by,

$$T_M = \frac{\int_{x_{L,M}}^{x_{M,N}} \rho c T(y) S dy}{\int_{x_{L,M}}^{x_{M,N}} \rho c S(y) dy}.$$

Given this, $x_M$ will be that location at which the isotherm with magnitude $T_M$ intersects the abcissa. Clearly, this position will be determined by the functional dependence of S on x. For any prescribed variation S(x), $x_M$ can be determined using (10), assuming known values of $T_{L,M}$ at $x_{L,M}$ and $T_{M,N}$ at $x_{M,N}$ (ref. 11).

Subject to these considerations, finite difference equation (6) for element M is to be replaced by,

$$\frac{T_L - T_M}{\int_{x_L}^{x_M} \frac{dy}{KA(y)}} - \frac{T_M - T_N}{\int_{x_M}^{x_N} \frac{dy}{KA(y)}} = \rho c \left( \int_{x_{L,M}}^{x_{M,N}} A(y) \, dy \right) \frac{\Delta T_M}{\Delta t},$$

where average temperature $\Delta T_M$ is associated with $x_M$ as described above.

Just as (7), this difference equation is sufficient in detail to obtain numerical solutions with prescribed boundary and initial conditions. The flow domain is visualized as a collection of flow tubes with curvilinear coordinate axes coinciding with suitably chosen flow lines.

For the special case of radial flow in a cylinder of thickness H for which $S(x) = 2\pi x H$, or for radial flow in a sphere for which $S(x) = 4\pi x^2$, (8) can be shown to lead to the appropriate differential equations by dividing through by volume, and assuming small spacing of elements.

**Discussion**

The continuum and the differential equation greatly influence modern scientific thinking.  There are two reasons.  Historically, when the absence of computing machines stood in the way of quantitatively solving practical problems of interest, the differential equation opened up a viable new way to make quantitative sense of the physical world.  Moreover, by casting problems algebraically, the differential equation enabled comprehension of whole classes of problems in a manner that would be impossible to achieve by solving problems numerically, on a case by case basis.  The remarkably rich body of knowledge that has accumulated over the past two centuries attests to the unique position that the differential equation and the continuum occupy in science.

Central to the continuum and the differential equation are the notions of points and gradients.  Both are abstract entities.  Physically, it is impossible to observe a point.  Nor is it possible to measure a gradient.  Physical measurements such as temperature are made with instruments of finite size and shape at discrete locations, separated by finite distances.  Gradients are inferred from these measurements.  The most powerful digital computer cannot truly represent an irrational number.  For these reasons, the difference equation provides a physically meaningful representation of natural phenomena, amenable to direct numerical estimations.  Thus, remarkably, the discrete world represented by the difference equation and the abstract world represented by the differential equation together help us comprehend the world around us in ways that neither can individually achieve.

The connections between the difference equation and the differential equation in probability and in heat diffusion presented above show that difference equations stemming from a discrete description of either phenomenon constitute the first stage in setting up the mathematical problem.  At this stage, one has the choice of solving a given problem numerically using a computing device.  Alternatively, one may invoke infinitesimals and continuous functions, and study the problem algebraically.  Each approach complements the other.

**Concluding Remark**

The differential equation and the continuum are so much part of our thinking that we start, *de facto*, with these notions even when confronted with discrete problems that have to be individually described and numerically solved.  We subordinate the discrete to the continuous.

The thoughts presented in this work should not be construed as a criticism or a negation of differential calculus. On the contrary, it is argued that the continuum and the differential equation should not negate the logical consistency of the discrete and the difference equation. The discrete and the continuous are two independent, complementary visualizations that the human mind has created to make sense of the complex physical world of our existence.

---

*T. N. Narasimhan, Professor Emeritus*

*Department of Materials Science and Engineering*

*210 Hearst Memorial Mining Building*

*University of California, Berkeley, Ca 94720-1760   USA*

*tnnarasimhan@LBL.gov*